\documentclass[10 pt,a4paper,twoside,reqno]{amsart}
\usepackage{amssymb, amsmath, amsthm, amsfonts, amscd, enumerate, verbatim, calc}

\newtheorem{theorem}{Theorem}[section]
\newtheorem{prop}[theorem]{Proposition}
\newtheorem{lemma}[theorem]{Lemma}
\newtheorem{cor}[theorem]{Corollary}
\theoremstyle{definition}
\newtheorem{defn}[theorem]{Definition}
\newtheorem{eg}[theorem]{Example}
\newtheorem{rem}[theorem]{Remark}

\title{$G_1$ class elements in a Banach algebra}
\author{S. H. Kulkarni}
\address{Indian Institute of Technology Palakkad, Ahalia Integrated Campus, Palakkad 678557. Kerala. (Formerly Department of Mathematics, Indian Institute of Technology - Madras, Chennai 600036)}
\email{shk@iitpkd.ac.in}
\begin{document}

\begin{abstract}
 Let $A$ be a complex unital  Banach algebra with unit $1$. An element  $a\in A$ is said  to be of \textit{$G_{1}$-class} if 
 $$\|(z-a)^{-1}\|=\frac{1}{\text{d}(z,\sigma(a))} \quad \forall z\in \mathbb{C}\setminus \sigma(a).$$
 Here $d(z, \sigma(a))$ denotes the distance between $z$ and the spectrum $\sigma(a)$ of $a$.
Some  examples of such elements are given and also some properties  are proved. It is shown that a 
$G_1$-class element is a scalar multiple of the unit $1$ if and only if its spectrum is a singleton set 
consisting of that scalar. It is proved that if $T$ is a $G_1$ class operator on a Banach space $X$, then every isolated point of $\sigma(T)$ is an eigenvalue of $T$. If, in addition, $\sigma(T)$ is finite, then $X$ is a direct sum of eigenspaces of $T$.

\end{abstract}

\maketitle
\section{Introduction}
Let $T$ be a normal operator on a complex Hilbert space $H$ and $\lambda$ a complex number  not lying in the spectrum $\sigma(T)$ of $T$. Then it is known that
the distance between $\lambda$ and $\sigma(T)$ is given by $\frac{1}{\|(\lambda I - T)^{-1}\|}$. It is also known that there are many other operators
that are not normal but still satisfy this property. Putnam called such operators as 
operators satisfying $G_1$ condition and investigated properties of such operators in 
\cite{putnam1}, \cite{putnam2}. In particular, he proved that if $T$ is a $G_1$ class operator, then every isolated point of $\sigma(T)$
is an eigenvalue of $T$ and every $G_1$ class operator on a finite dimensional Hilbert space is normal.

In this note we extend this concept of $G_1$ class operators to operators on a Banach space and more
generally to elements of a complex Banach algebra and investigate the properties of such elements.
The next section contains some preliminary definitions and results that are used throughout. In 
Section 3, we give definition of a $G_1$ class element in a complex unital Banach algebra, give some examples and prove a few elementary properties of such elements. In particular, it is proved that every
element of a uniform algebra is of $G_1$ class and conversely if every element of a complex unital Banach algebra  $A$ is of $G_{1}$ class, then $A$ is commutative, semisimple and hence isomorphic and homeomorphic to a uniform algebra. The last section deals with the spectral properties of $G_1$ class elements and contains the main results of this note. In particular, it is proved that if $T$ is a $G_1$ class operator on a Banach space $X$, then every isolated point of $\sigma(T)$ is an eigenvalue of $T$. Further, if, in 
addition, $\sigma(T)$ is finite, then $X$ is a direct sum of eigenspaces of $T$. In this sense $T$ is ``diagonalizable" and hence this result can be considered to be an analogue of the Spectral Theorem for such operators. 

An overall aim of such a study can be to obtain an analogue of the Spectral Theorem for $G_{1}$ class operators. Though at present we are far away from this goal, the present results can be considered a small step in that direction. Next natural step should be to try to prove a similar result for compact 
operators of $G_{1}$ class. Another way of looking at this study is an attempt to answer the following question:``To what extent does the spectrum of an element determine the element?" This question has a long and interesting history. It has appeared under different names at different times such as ``Spectral characterizations", ``hearing the shape of a drum",\cite{kac} ``$T = I$ problem" \cite{Zang} etc. The results in this note say that the spectrum of a $G_1$ class element gives a fairly good information about   
that element.

We shall use the following notations throughout this article.
 Let\\
  $B(w, r):=\{z\in\mathbb{C}:|z- w|< r\}$, the open disc with the centre at $w$ and radius $r$, \\ 
  ${D(z_{0},r):=\{z\in\mathbb{C}:|z-z_{0}|\leq r\}}$, the closed disc with the centre at $w$ and radius $r$, \\
 ${A+D(0,r)= \bigcup \limits_{a\in A}D(a;r)\text{ for } A\subseteq \mathbb{C}}$
and ${d(z,K)=\text{inf}\{|z-k|:k\in K\}}$, the distance between a complex number $z$  and a  closed set $K\subseteq \mathbb{C}$.\\
Let $\delta \Omega$ denote the boundary of a set $\Omega\subseteq \mathbb{C}$.\\
$\mathbb{C}^{n\times n}$ denotes the space of square matrices of order $n$ and $B(X)$ denotes the set of bounded linear operators on a Banach space $X$.

\section{Preliminaries}
Since our main objects of study are  certain elements in a Banach algebra, we shall review some definitions related to a Banach algebra. Many of these definitions can be found in the book \cite{bonsall}. Some material in this section is also available in the review article \cite{shk}.

\begin{defn}{\bf Spectrum}:
Let $A$ be a complex unital Banach algebra with unit $1$. For $\lambda \in \mathbb{C}, \, \lambda.1$ is identified with $\lambda$.
Let $\text{Inv}(A)=\{x\in A:x \text{ is invertible in }A\}$ and $\text{Sing}(A)=\{x\in A:x \text{ is not invertible in }A\}.$
The {\it spectrum} of an element $a\in A$ is defined as:
\begin{displaymath}\sigma (a):=\{ \lambda \in \mathbb{C}:\lambda-a \in \text{Sing}(A)\} \end{displaymath}
The {\it spectral radius} of an element $a$ is defined as:
\begin{displaymath}r(a):=\text{sup} \{ |\lambda|:\lambda \in \sigma(a)\} \end{displaymath}

Its value is also given by  the Spectral Radius Formula,

$$ r(a) = \lim_{n \rightarrow \infty} \|a^n\|^{\frac{1}{n}} = \inf_n  \|a^n\|^{\frac{1}{n}}  $$
The complement of the spectrum of an element $a$ is called the {\it resolvent set of $a$} and is denoted by $\rho(a)$.
\end{defn}
Thus when $A = C(X)$, the algebra of all continuous complex valued functions on a compact Hausdorff space $X$ and $f \in A$,then the spectrum  $\sigma(f)$ of $f$ coincides with the range of $f$.

Similarly when $A = \mathbb{C}^{n \times n}$, the algebra of all square matrices of order $n$ with complex entries and $M \in A$, the spectrum $\sigma(M)$ of $M$ is the set of all eigenvalues of $M$.

\begin{defn}{\bf Numerical Range}
Let $A$ be a Banach algebra and $a \in A$. The {\it numerical range}  of $a$ is defined by \begin{displaymath}V(a):=\{f(a):f\in A',f(1)=1=\|f\| \},\end{displaymath} where $A'$ denotes the dual space of $A$, the space of all continuous linear functionals on $A$..

The {\it numerical radius} $\nu(a)$ is defined as \\
$$ \nu(a) := \sup \{|\lambda|: \lambda \in V(a)\}    $$

 Let $A$ be a Banach algebra and $a\in A$. Then $a$ is said to be {\it Hermitian} if $V(a)\subseteq\mathbb{R}$.
\end{defn}
If $A$ is a is a $C^*$ algebra(also known as $B^*$ algebra), then an element $a\in A$ is Hermitian if and only if 
it is self-adjoint. \cite{bonsall}

\begin{defn}{\bf Spatial Numerical Range} 
 
 Let $X$ be a Banach space and $T\in B(X)$. Let $X^{'}$ denote  the dual space of $X$. The {\it spatial numerical range} of $T$ is defined by 
 \begin{displaymath}
  W(T)=\{f(Tx):f\in X^{'}, \|f\|=f(x)=1= \|x\| \}.
 \end{displaymath}
\end{defn}

For an operator $T$ on a Banach space $X$, the spatial numerical range $W(T)$ and the numerical range $V(T)$, where $T$ is regarded as 
an element of the Banach algebra $B(X)$, are related by the following: 
\begin{displaymath}
 \overline {\text{Co}}\,W(T)=V(T)
\end{displaymath}
where $\overline {\text{Co}}\,E$ denotes the closure of the  convex hull of $E\subseteq\mathbb{C}$. 

The following theorem gives the relation between the spectrum and numerical range.

\begin{theorem}
 Let $A$ be a complex unital Banach algebra with unit $1$ and $a \in A$. \\
 Then the numerical range $V(a)$ is a closed convex set containing $\sigma(a)$. Thus\\
 $\overline{Co}(\sigma(a)) \subseteq V(a)$. Hence\\
 $r(a) \leq \nu(a) \leq \|a\| \leq e \nu(a).  $
 \end{theorem}
A proof of this can be found in \cite{bonsall}.

\begin{cor}
Let $A$ be a complex unital Banach algebra with unit $1$ and $a \in A$.
If $a$ is Hermitian, then $\sigma(a) \subseteq \mathbb{R}$.
\end{cor}

We now discuss another important and popular set related to the spectrum,
namely pseudospectrum. We begin with its definition.

\begin{defn}{\bf Pseudospectrum} 
Let $A$ be a complex Banach algebra, $a\in A$ and $\epsilon>0$.
The $\epsilon$-pseudospectrum $\Lambda_{\epsilon} (a)$ of $a$ is defined by \begin{displaymath}{\Lambda_{\epsilon} (a):=\{ \lambda \in \mathbb{C}:\|(\lambda-a)^{-1}\| \geq \epsilon^{-1}\}}\end{displaymath}
with the convention that $\|(\lambda-a)^{-1}\|=\infty$ if $\lambda-a$ is not invertible.
\end{defn}

This definition and many results in this section can be found in \cite{arundhathi1}. The book \cite{trefethen} is 
a standard reference on Pseudospectrum. It contains a good amount of information about the idea of pseudospectrum,
(especially in the context of matrices and operators), historical remarks and applications to various fields. Another useful 
source is the website \cite{pseudoweb}.

The following theorems establish the relationships between the spectrum, the $\epsilon$-pseudospectrum
and the numerical range of an element of a Banach algebra.
\begin{theorem} \label{thm:relationship}
Let $A$ be a Banach algebra, $a\in A$ and $\epsilon>0$. Then
\begin{equation}
d(\lambda,V(a))\leq\frac{1}{\|(\lambda-a)^{-1}\|}\leq d(\lambda, \sigma(a)) \quad  \forall \lambda\in \mathbb{C}\setminus\sigma(a).
\end{equation}
Thus
\begin{equation}\sigma(a)+D(0;\epsilon)\subseteq \Lambda_{\epsilon}(a)\subseteq V(a)+D(0;\epsilon).\end{equation} 
\end{theorem}
A proof of this Theorem can be found in \cite{arundhathi1}.

The following theorem gives the basic information  about  the analytical functional calculus for elements of a Banach algebra. 
\begin{theorem}\label{Spectral mapping theorem}
Let $A$ be a Banach algebra and $a\in A$. Let $\Omega\subseteq \mathbb{C}$ be an open neighbourhood of $\sigma(a)$ and $\Gamma$
be a contour that surrounds $\sigma(a)$ in $\Omega$. Let $H(\Omega)$ denote the set of all  analytic functions in $\Omega$ and let $P(\Omega)$ denote the set of all polynomials in $z$ with $z \in \Omega$. We recall the definition of $\tilde{f}(a)$ in the analytical functional calculus as 
 \begin{equation}\tilde{f}(a)=\frac{1}{2\pi i} \int\limits_{\Gamma}(z-a)^{-1}f(z)dz \end{equation}
 Then the map $f \rightarrow \tilde{f}(a)$ is a homomorphism from $H(\Omega)$ into $A$ that extends the natural homomorphism $p \rightarrow p(a)$ of $P(\Omega)$ into $A$ and 
 $$\sigma(\tilde{f}(a)) = \{f(z): z \in \sigma(a)\}$$
 \end{theorem}
 A proof of this Theorem can be found in \cite{bonsall}.

\section{$G_1$ class elements}

In this section, we give definition, some examples and elementary properties of 
$G_1$ class elements. It is possible to view this definition as motivated by considering 
 the question of equality in some of the inclusions given in Theorem \ref{thm:relationship}. 
\begin{defn}
 Let $A$ be a Banach algebra and $a\in A$. We define $a$ to be of \textit{$G_{1}$-class} if 
 \begin{equation}\|(z-a)^{-1}\|=\frac{1}{\text{d}(z,\sigma(a))} \quad \forall z\in \mathbb{C}\setminus \sigma(a).\end{equation}
\end{defn}

 \begin{rem}
 
 The idea of $G_{1}$-class was introduced by Putnam who defined it for operators on Hilbert spaces.
 (See \cite{putnam1},\cite{putnam2}.) 
 It is known that the $G_{1}$-class properly contains the class of seminormal operators
 (that is, the operators satisfying $TT^{*}\leq T^{*}T$ or $T^{*}T\leq TT^{*}$) and this class properly contains the class of normal operators.
 Using the Gelfand- Naimark theorem \cite{bonsall}, we can make similar statements about elements in a $C^{*}$ algebra. \\
 $G_{1}$-class operators on a finite dimensional Hilbert space are normal\cite{putnam1}.\\
 
 In particular, normal elements are hyponormal. In general, the equation (4) may hold, for every $z\in \mathbb{C}\setminus \sigma(a)$,
 for an element $a$ of a $C^{\ast}$-algebra even though $a$ is not normal. \\
 For example, we may consider the right shift operator $R$ on $\ell ^{2}(\mathbb{N})$. 
 It is not normal but $\Lambda_{\epsilon}(R)=\sigma(R)+D(0;\epsilon)=D(0;1+\epsilon) \, \forall \epsilon>0$.
The operator $R$ is, however, a hyponormal operator.
\end{rem}

We now deal with a natural question: What are $G_1$ class elements in an arbitrary Banach algebra?

The following lemma is elementary and gives a characterization of a $G_1$ class element in terms of its pseudospectrum. 
\begin{lemma}\label{lem:equal}
Let $A$ be a Banach algebra and $a\in A$. Then \begin{equation}\label{eq:equal}\Lambda_{\epsilon}(a)=\sigma(a)+D(0;\epsilon) \quad \forall \epsilon>0\end{equation}
iff $a$ is of $G_{1}$-class. 
\end{lemma}
A proof of this Lemma can be found in \cite{arundhathi1}.

As one may expect, most natural candidates to be $G_1$ class elements are scalars, that is, scalar multiples of the identity $1$. 

\begin{theorem}\label{thm:scalar}
Let $A$ be a complex Banach algebra with unit $1$ and $a\in A$. \\
(i) If $a = \mu $ for some complex number $\mu$, then $a$ is of $G_1$ 
class and $\sigma(a) = \{\mu\}$.\\
(ii) If  $a$ is of $G_1$ class, then $\alpha a + \beta$ is also of $G_1$ 
class for every complex numbers $\alpha, \beta$.\\
(iii) If $a$ is of $G_1$ class and $\sigma(a) = \{\mu\}$, then $a = \mu $.
\end{theorem}

A proof of this is straight forward. It also follows easily from \ref{lem:equal} and Corollary 3.17 of \cite{arundhathi1}. We include it here for the sake of completeness.
\begin{proof}
(i) Let $a = \mu $ for some complex number $\mu$. Then clearly $\sigma(a) = \{\mu\}$. Hence for all $z\in \mathbb{C}\setminus \sigma(a)$, we have $z \neq \mu$. Thus $\|(z-a)^{-1}\|=\frac{1}{|z -\mu|} = \frac{1}{\text{d}(z,\sigma(a))}$. This shows that  $a$ is of $G_1$ class.\\
(ii) Next suppose that  $a$ is of $G_1$ class and $b = \alpha a + \beta$
for some complex numbers $\alpha, \beta$. We want to prove that $b$ is of $G_1$ class. If $\alpha = 0$, then it follows from (i). So assume that $\alpha \neq 0$. Let $w \notin \sigma(b) = \{\alpha z + \beta : z \in \sigma(a)\}$. Then $z:= \frac{w - \beta}{\alpha} \notin \sigma(a)$ and since 
$a$ is of $G_1$ class, $\|(z-a)^{-1}\|=\frac{1}{\text{d}(z,\sigma(a))}$.
Now $\|(w -b)^{-1}\| = \|(\alpha z + \beta - (\alpha a + \beta))^{-1}\|
= \frac{1}{|\alpha|} \|(z-a)^{-1}\|= \frac{1}{|\alpha|d(z,\sigma(a)}= 
\frac{1}{d(\alpha z, \sigma(\alpha a))} = \frac{1}{d(w, \sigma(b))}$. This
shows that $b$ is of $G_1$ class. \\
(iii) Suppose $a$ is of $G_1$ class and $\sigma(a) = \{\mu\}$. Let $b = a - \mu$. Then by (ii), $b$ is of $G_1$ class and $\sigma(b) = \{0\}$. Let $\epsilon > 0$ and $C$ denote the circle with the centre at $0$ and radius $\epsilon$ traced anticlockwise. Then for every $z \in C$, $\|(z - b)^{-1}\| = \frac{1}{d(z, \sigma(b)} = \frac{1}{|z-0|}= \frac{1}{\epsilon}$. Also\\
$$b = \frac{1}{2\pi i} \int_C z (z - b)^{-1} dz$$
Hence $\|b\| \leq \frac{1}{2\pi}2\pi\epsilon \epsilon \frac{1}{\epsilon} = \epsilon$. Since this holds for every $\epsilon >0$, we have $b = 0$, that is $a = \mu$.
\end{proof}

\begin{rem}
The above Theorem has a relevance in the context of a
  very well known classical problem in operator theory known as {\it ``$T=I$? problem"}. This problem  asks the following question: {\it Let $T$ be an operator on a Banach space. Suppose $\sigma(T)=\{1\}$. Under what additional conditions can we conclude $T=I$?} A survey article \cite{Zang} contains details of many classical results about this problem.

 From the above Theorem it follows that if $T$ is of $G_1$ class and  
 $\sigma(T)=\{1\}$,  then
 we can conclude that   $T=I$. In other words ``$T$ is of $G_1$ class'' works as an additional condition in the $``T=I$ problem".
 \end{rem}

Next we show  that every Hermitian idempotent element is of $G_1$ class. A version of this result 
was included in the thesis \cite{arundhathithesis}.

\begin{theorem}\label{thm:hermitianidempotent}
Let $A$ be a complex unital Banach algebra with unit $1$ and $a\in A$. If $a$ is 
a Hermitian idempotent element, then $a$  is of $G_1$ class. Also, if 
$a$  is of $G_1$ class and $\sigma(a) \subseteq \{0, 1\}$, then $a$ is a Hermitian idempotent. 

\end{theorem}
\begin{proof}
Supose $a$ is a Hermitian idempotent element. If $a = 0$ or $a = 1$, then 
 $a$  is of $G_1$ class by (i) of Theorem \ref{thm:scalar}. Next, let $a \neq 0, 1$. Then $\sigma(a) = \{0, 1\}$ and by Theorem 1.10.17 of \cite{bonsall}, $\|a\| = r(a) = 1$. Now Corollary 3.18 of \cite{arundhathi1}
 implies that $\Lambda_{\epsilon}(a) = D(0,\epsilon) \cup D(1, \epsilon)$ for every $\epsilon > 0$. Hence $a$ is of $G_1$ class by Lemma \ref{lem:equal}.
 
 Next suppose $a$  is of $G_1$ class and $\sigma(a) \subseteq \{0, 1\}$.
 If $\sigma(a) = \{0\}$, then $a = 0$ by (ii) of Theorem \ref{thm:scalar}.
 Similarly, if $\sigma(a) = \{1\}$, then $a = 1$. So assume that  $\sigma(a) =  \{0, 1\}$. Then by Lemma \ref{lem:equal},  $\Lambda_{\epsilon}(a) = D(0,\epsilon) \cup D(1, \epsilon)$ for every $\epsilon > 0$. Hence by  3.18 of \cite{arundhathi1}, $a$ is a Hermitian idempotent element.
\end{proof}

The abundance or scarcity of $G_1$ class elements in a given Banach algebra depends on the nature of that Banach algebra. There exist extreme cases, that is, there are Banach algebras in which every element is of $G_1$ class.
On the other hand, there are also Banach algebras in which the scalars are the only elements of $G_1$ class. We shall see examples of both types below.
Before that, we need to review a relation between the spectrum and numerical range of an element of $G_1$ class. Recall that the numerical range of an element of a Banach algebra is a compact convex subset of $\mathbb{C}$ containing its spectrum,
and hence it also contains the closure of the convex hull of the spectrum. The next proposition shows that  the equality holds in case of elements of $G_1$ class.\\
 \begin{prop}\label{prop:convexhull}
  Let $A$ be a complex unital Banach algebra and $a\in A$. Suppose $a$ is of $G_{1}$-class.
  Then $V(a)=\overline{\text{Co }}(\sigma(a))$, the closure of the convex hull of the spectrum of $a$  and $\|a\|\leq e \,r(a)$.
   \end{prop}
   A proof of this can be found in \cite{arundhathi1}.
   
   \begin{cor}
   Let $A$ be a complex unital Banach algebra. Suppose $a \in A$ is of $G_{1}$-class and $\sigma(a) \subseteq \mathbb{R}$. Then $a$ is Hermitian.
   \end{cor}
   
 It is shown in the next theorem that every element in a uniform algebra is of $G_{1}$ class. Also a  partial converse of this statement is proved. We may recall that a {\it uniform algebra} is a unital Banach algera satisfying $\|a\|^2 = \|a^2\|$ for every $a \in A$. Every complex uniform algebra is commutative by a theorem of Hirschfeld and Zelazko \cite{bonsall}. Then it follows by Gelfand theory \cite{bonsall} that such an algebra is isomertically isomorphic to a {\it function algebra}, that is, a uniformly closed subalgebra of $C(X)$ that contains the constant function $1$ and separates the points of $X$, where $X$ is the maximal ideal space of $A$.

\begin{theorem}\label{thm:G_1 uniform} (See also Theorem 3.15 of \cite{arundhathi1})
 Let $A$ be a complex unital Banach algebra with unit $1$. \\
 (i) If $A$ is a uniform algebra, then every element in $A$ is of $G_{1}$class.\\
 (ii) If every element of $A$ is of $G_{1}$ class,
 then $A$ is commutative, semisimple and hence isomorphic and homeomorphic to a uniform algebra.
\end{theorem}
\begin{proof}
(i) The Spectral Radius Formula implies that $\|a\| = r(a)$ for every $a \in A$. Now let $a \in A$ and $\lambda \notin \sigma(a)$. Then \\
\begin{align*}
\|(\lambda - a)^{-1}\| &= r((\lambda - a)^{-1})\\
 &= \sup\{|z|: z \in \sigma((\lambda - a)^{-1}) \} \\
  &= \sup\{\frac{1}{|\lambda - \mu|}:\mu \in \sigma(a)\}\\
   &= \frac{1}{\inf\{|\lambda - \mu|:\mu \in \sigma(a)\}}\\
    &= \frac{1}{d(\lambda, \sigma(a)} \end{align*}.
This shows that $a$ is of $G_1$ class.\\

 (ii) By  Proposition \ref{prop:convexhull}, $\|a\|\leq e r(a)$ for all $a\in A$.
Hence $A$ is commutative by a theorem of Hirschfeld and Zelazko \cite{bonsall}. Also, the condition 
$\|a\|\leq e r(a)$ for all $a\in A$ implies that $A$ is semisimple and hence the spectral radius $r(.)$ is a norm on $A$. Clearly, $r(a^2) = (r(a))^2$ for every $a \in A$. Hence  $A$ is a uniform algebra under this norm. Also the inequality   $r(a) \leq \|a\|\leq e r(a)$ for all $a\in A$ implies that
the identity map is a homeomorphism between these two algebras.
\end{proof}

Next we consider an example of a Banach algebra in which scalars are the only elements of $G_1$ class.

\begin{eg}(See also Example 2.16 and Remark 2.20 of \cite{kousik1})

Let $A = \{ a \in \mathbb{C}^{2 \times 2}: a = \begin{bmatrix}
\alpha & \beta\\ 0 & \alpha
\end{bmatrix} \} $
with the norm given by $\|a\| = |\alpha| + \|\beta|$.
Suppose $ a = \begin{bmatrix}
\alpha & \beta\\ 0 & \alpha
\end{bmatrix} \in A$ 
is of $G_1$ class. Then since $\sigma(a) = \{\alpha\}$, it follows by Theorem \ref{thm:scalar}(iii) that $a = \alpha$. (This means $\beta = 0$.)

\end{eg}

\section{Spectral properties of $G_1$ class elements}
In this section, we show that $G_1$ class elements have some properties that are very similar to the properties of normal operators on a complex Hilbert space. For example, if $H$ is a complex Hilbert space, $T$ is 
a normal operator on $ H$ and $\lambda$ is an isolated point of $\sigma(T)$, then $\lambda$ is an eigenvalue of $T$. We show that a similar property holds for a bounded operator of $G_1$ class on a Banach space. For that we need the following theorem about isolated points 
of the spectrum of a $G_1$ class element in a Banach algebra. 

\begin{theorem}\label{thm:isolated}
Let $A$ be a complex unital Banach algebra with unit $1$. Suppose $a$ is of $G_{1}$-class and $\lambda$ is an isolated point of $\sigma(a)$. Then there exists an idempotent element $ e \in A$ such that $ae = \lambda e$ and 
$\|e\| = 1$. 
\end{theorem}

\begin{proof}
If $\sigma(a) = \{\lambda\}$, then by \ref{thm:scalar}(iii), $a = \lambda$ and we can take $e = 1$. 

Next assume that $\sigma(a) \setminus \{\lambda\}$ is nonempty. Let $D_1$ and $D_2$ be disjoint open neighbourhoods of $\lambda$ and $\sigma(a) \setminus \{\lambda\}$ respectively. Define 
$$f(z) = \begin{cases} 1 \quad \text{if}\quad z \in D_1\\ 0 \quad \text{if} \quad z \in D_2 \end{cases}    $$
Then $f$ is analytic in $D_1 \cup D_2$. Let $e = \tilde{f}(a)$. Then since $f^2 = f$, we have $e^2 = e$, that is, $e$ is an idempotent element and $\|e\| \geq 1$. To prove other assertions, choose $\epsilon > 0$ in such a way that for every 
$z \in \Gamma_1 := \{w \in \mathbb{C}: |w - \lambda| = \epsilon\}$, $\lambda$ is the nearest point of  $\sigma(a)$ and $\Gamma_1 \subseteq D_1$. Then for every such $z$, $d(z, \sigma(a)) = |z - \lambda| = \epsilon$, hence $\|(z -a)^{-1}\| = \frac{1}{\epsilon}$. Now let $\Gamma_2$ be any closed curve lying in $D_2$ and enclosing $\sigma(a) \setminus \{\lambda\}$ and let $\Gamma = \Gamma_1 \cup \Gamma_2$. Then 
$$e = \tilde{f}(a) = \frac{1}{2\pi i}\int_{\Gamma}f(z)(z-a)^{-1}dz = \frac{1}{2\pi i}\int_{\Gamma_1}(z-a)^{-1}dz$$
Hence 
$$ \|e\| \leq \frac{1}{2\pi} \frac{1}{\epsilon}2\pi\epsilon =1 $$
This shows that $\|e\| = 1$.

Now define $g(z) = (z - \lambda)f(z)$. Then $|g(z)| \leq \epsilon$ for all
$z \in \Gamma_1$. Note that
$$ ae - \lambda e = \tilde{g}(a) = \frac{1}{2\pi i}\int_{\Gamma}g(z)(z-a)^{-1}dz = \frac{1}{2\pi i}\int_{\Gamma_1}g(z)(z-a)^{-1}dz$$   
Hence 
$$ \|ae - \lambda e\| \leq \frac{1}{2\pi} \, \epsilon \, \frac{1}{\epsilon}\, 2\pi\epsilon =\epsilon $$
Since this holds for every $\epsilon > 0$, we have $ae - \lambda e = 0$.

\end{proof}

\begin{cor}\label{cor:eigenvalue}
Let $X$ be a complex Banach space, $T \in B(X)$ be of $G_1$ class and $\lambda$ be an isolated point of $\sigma(T)$. Then $\lambda$ is an eigenvalue of $T$. 
\end{cor}

\begin{proof}
By Theorem \ref{thm:isolated}, there exists an idempotent element $P \in B(X)$ such that $\|P\| = 1$ and $TP = \lambda P$. Clearly $P$ is a nonzero 
projection operator on $X$. Let $x \neq 0$ be an element of the range $R(P)$ 
of $P$. Then $P(x) = x$. Hence $T(x) = TP(x) = \lambda P(x) = \lambda x $.
Thus $\lambda$ is an eigenvalue of $T$.
\end{proof}

Some ideas in the proof of the next theorem can be compared with the proof of Theorem C
in \cite{stampfli} that deals with similar results about hyponormal operators on a Hilbert space.
\begin{theorem}\label{thm:finite spectrum}
Let $A$ be a complex unital Banach algebra with unit $1$. Suppose $a$ is of $G_{1}$-class and $\sigma(a) = \{\lambda_1, \ldots, \lambda_m\}$ is finite. 
Then there exist idempotent elements $e_1, \ldots, e_m $ such that
\begin{enumerate}

\item  $\|e_j\|  = 1$, $ae_j = \lambda_je_j$ for $j = 1, \ldots, m$,  $e_j e_k = 0$ for $j \neq k$,
$$ e_1 + \ldots +  e_m = 1 $$
and  $$ a = \lambda_1 e_1 + \ldots + \lambda_m e_m.  $$
\item If $p$ is any polynomial, then
$$ p(a) = p( \lambda_1) e_1 + \ldots + p(\lambda_m) e_m.  $$
\item In particular, $$(a - \lambda_1) \ldots (a - \lambda_m) = 0.$$
\item If $\lambda$ is a complex number such that $\lambda \neq \lambda_j$ 
for $j = 1, \ldots, m$, then 
$$(\lambda - a)^{-1} = \frac{1}{\lambda - \lambda_1}e_1 + \ldots  + \frac{1}{\lambda - \lambda_m}e_m. $$
\item If a function $f$ is analytic in a neighbourhood of $\sigma(a)$, then
$$\tilde{f}(a) = f(\lambda_1)e_1 + \ldots +   f(\lambda_m)e_m  $$
\end{enumerate}
\end{theorem}
\begin{proof}
If $m = 1$, then by Theorem \ref{thm:scalar}(iii), $a = \lambda_1$. Hence we 
can take $e_1 = 1$ and all the conclusions follow trivially.
Next we assume $m > 1$. Let $D_1, \ldots, D_m$ be mutually disjoint neighbourhoods of $\lambda_1, \ldots, \lambda_m\ $ respectively and let
$D = \cup_{j=1}^{m} D_j$. Now for each $j = 1, \ldots, m$, define a function $f_j$ on  $D$ by
$$ f_j(z) = \begin{cases} 1 \quad \mbox{if} \quad z \in D_j\\
0 \quad \mbox{if} \quad z \notin D_j \end{cases}  $$
Let $e_j = \tilde{f_j}(a) $. Then it follows as in Theorem \ref{thm:isolated}
that each $e_j$ is an idempotent, $\|e_j\| = 1$ and $ae_j = \lambda_j e_j$.
Since for $j \neq k$, $f_j f_k = 0$, we have $e_j e_k = 0$. \\
Further $f_1 + \ldots + f_m = 1$ implies $ e_1 + \ldots +  e_m = 1 $.\\
Next
\begin{align*}
a &= a 1 \\ &= a ( e_1 + \ldots +  e_m)\\ &= a e_1 + \ldots + a e_m  \\
&= \lambda_1 e_1 + \ldots + \lambda_m e_m.
\end{align*} 
This proves (1).\\
Now since $e_j^2 = e_j$ for each $j$ and $e_j e_k = 0$ for $j \neq k$, we have
$$ a^2 = \lambda_1^2 e_1 + \ldots + \lambda_m^2 e_m   $$

and in general for any power $k$,
$$ a^k = \lambda_1^k e_1 + \ldots + \lambda_m^k e_m.   $$
It follows easily from this that for any polynomial $p$, we have 
$$ p(a) = p( \lambda_1) e_1 + \ldots + p(\lambda_m) e_m.  $$
Thus (2) is proved.\\
Now consider the polynomial $p$ given by  $p(z) = (z - \lambda_1) \ldots (z - \lambda_m) $. Then $p(\lambda_j) = 0$ for each $j$. Hence $p(a) = 0$, that is, $(a - \lambda_1) \ldots (a - \lambda_m) = 0$.
This completes the proof of (3).\\
Now suppose $\lambda$ is a complex number such that $\lambda \neq \lambda_j$ 
for $j = 1, \ldots, m$. Let 
$$b = \frac{1}{\lambda - \lambda_1}e_1 + \ldots  + \frac{1}{\lambda - \lambda_m}e_m. $$
Then in view of (1), we have
\begin{align*}
(\lambda - a) b &= [(\lambda - \lambda_1)e_1 + \ldots + (\lambda - \lambda_m)e_m][\frac{1}{\lambda - \lambda_1}e_1 + \ldots  + \frac{1}{\lambda - \lambda_m}e_m]\\
&= 1
\end{align*}

Similarly, we can prove $b(\lambda - a)= 1$ implying (4).\\
Next suppose a function $f$ is analytic in a neighbourhood  $\Omega$ of  $\sigma(a)$ and 
$\Gamma$ is a closed curve lying in  $\Omega$ and surrounding   $\sigma(a)$. Then

\begin{align*}
\tilde{f}(a) &= \frac{1}{2\pi i}\int_{\Gamma}f(z)(z - a)^{-1} dz \\
&= \frac{1}{2\pi i}\int_{\Gamma}f(z)[\frac{1}{z - \lambda_1}e_1 + \ldots  + \frac{1}{z - \lambda_m}e_m] dz\\
&= (\frac{1}{2\pi i}\int_{\Gamma}\frac{f(z)}{z - \lambda_1}dz)e_1 + \ldots + (\frac{1}{2\pi i}\int_{\Gamma}\frac{f(z)}{z - \lambda_m}dz)e_m\\
&= f(\lambda_1)e_1 + \ldots +   f(\lambda_m)e_m
\end{align*}
\end{proof}

\begin{rem}
Note that the conclusions (2) and (4) of the above Theorem are special cases of (5).
\end{rem}

Now we apply the above Theorem to a bounded operator on a Banach space.

\begin{theorem}\label{thm:spectral thm}
Let $X$ be a complex Banach space. Suppose  $T \in B(X)$ is of $G_1$ class and  $\sigma(T) = \{\lambda_1, \ldots, \lambda_m\}$ is finite. Then

\begin{enumerate}
\item  Each $\lambda_j$ is an eigenvalue of $T$. In fact, there exist projections $P_j$
such that for each $j$, the range of $P_j$ is the  eigenspace corresponding to the eigenvalue $\lambda_j$ and $X$ is the direct sum of these eigenspaces. In other words, $T$ is ``{\it diagonalizable}". Also $\|P_j\| = 1$ and $TP_j = \lambda_j P_j$ for each $j$, $P_j P_k = 0$ for $j \neq k$, 

$$ P_1 + \ldots +  P_m = I $$
and  $$ T = \lambda_1 P_1 + \ldots + \lambda_m P_m.  $$
\item  $$(T - \lambda_1I) \ldots (T - \lambda_mI) = 0.$$
\item If a function $f$ is analytic in a neighbourhood of $\sigma(T)$, then
$$\tilde{f}(T) = f(\lambda_1)P_1 + \ldots +   f(\lambda_m)P_m  $$

\end{enumerate}
\end{theorem}
\begin{proof}
It follows from Corollary \ref{cor:eigenvalue} that each $\lambda_j$ is an eigenvalue of $T$.
The existence and properties of projections $P_j$ follow from Theorem \ref{thm:finite spectrum}.
Let $X_j =R(P_j) $, the range of $P_j$. The property $TP_j = \lambda_j P_j$ implies that 
$X_j$ is the eigenspace of $T$ corresponding to the eigenvalue $\lambda_j$ for each $j$. Also $P_j P_k = 0$ for $j \neq k$
implies that $X_j\cap X_k = \{0\} $ for $j \neq k$. It follows from
$$ P_1 + \ldots +  P_m = I $$
that $X$ is the sum of $X_j$. This shows that $X$ is the direct sum of these eigenspaces.

\end{proof}

\begin{rem}
Let $X$ and $T$ be as in the above Theorem. Since the conclusion (1) says that $X$ has a basis consisting 
of eigenvectors of $T$ and $T$ is a linear combination of projections, it can be called Spectral Theorem for such operators. Similarly, the conclusion (2) is an analogue of the Caley-Hamilton Theorem. If, in particular, $X$ is a Hilbert space, then every projection of norm $1$ is orthogonal and hence Hermitian(self-adjoint). Thus each $P_j$ is self-adjoint and hence $T$ is normal. This result is also proved in \cite{putnam2}. 

Suppose $X$ is finite dimensional. Then the above Theorem says that every $G_1$ class operator on $X$
is diagonalizable.
\end{rem}

\centerline{----------------------------------------------------------}
\end{document}